\documentclass[11pt]{article}
\usepackage{amsmath,amssymb,mathrsfs}
\usepackage{amsfonts}
\usepackage{sectsty}
\usepackage{hyperref}
\input xy
\xyoption{all}

\subsectionfont{\mdseries\itshape}

\newtheorem{Thm}{Theorem}[section]
\newtheorem{Prop}[Thm]{Proposition}
\newtheorem{Lemma}[Thm]{Lemma}
\newtheorem{Cor}[Thm]{Corollary}

\newcommand{\pf}{\noindent{\em Proof.}\ }
\newcommand{\qed}{\hfill $\Box$}

\renewcommand{\ni}{\noindent}

\newcommand{\ZZ}{\mathbb{Z}}
\newcommand{\QQ}{\mathbb{Q}}
\newcommand{\CC}{\mathbb{C}}
\newcommand{\FF}{\mathbb{F}}
\newcommand{\PP}{\mathbb{P}}
\renewcommand{\AA}{\mathbb{A}}
\newcommand{\tensor}{\otimes}
\newcommand{\uc}{\cup}
\newcommand{\dc}{\cap}
\newcommand{\rc}{\subset}
\newcommand{\lc}{\supset}

\newcommand{\Spec}{\operatorname{Spec}}

\newcommand{\tha}{\theta}
\newcommand{\om}{\omega}
\newcommand{\Fil}{\operatorname{Fil}}
\newcommand{\Gr}{\operatorname{Gr}}
\newcommand{\Has}{{\mathfrak H}}
\newcommand{\HH}{\mathcal{H}}
\newcommand{\MM}{{\mathcal M}}
\newcommand{\LL}{\mathcal{L}}
\newcommand{\JJ}{\mathcal{J}}
\newcommand{\NN}{\mathcal{N}}
\newcommand{\OO}{\mathcal{O}}

\newcommand{\lam}{\lambda}
\newcommand{\Lam}{\Lambda}
\newcommand{\Hom}{\operatorname{Hom}}
\newcommand{\ch}{\check}
\newcommand{\llsb}{[\hspace{-1.5pt}[}
\newcommand{\rrsb}{]\hspace{-1.5pt}]}
\newcommand{\lan}{\langle}
\newcommand{\ran}{\rangle}

\begin{document}

\title{Notes on Calabi-Yau Ordinary Differential Equations\footnote{
This work is supported by a grant
from the National Science Council of Taiwan.}}
\author{\sc Jeng-Daw Yu}
\date{2008.10.22}
\maketitle

\begin{abstract}
We investigate
the structures of Calabi-Yau differential equations
and the relations to the arithmetic of the pencils
of Calabi-Yau varieties behind the equations.
This provides explanations of some observations and computations
in the recent paper \cite{SvS}.
\end{abstract}

\setcounter{section}{-1}
\section{Introduction}

This note may be regarded as an appendix to the paper
\cite{SvS} of Samol and van Straten.
In that paper,
the authors study the variation of the frobenius action
on the cohomology of a Calabi-Yau variety
along a projective line (over a finite field)
via the expansion at a totally degenerate point
of a certain period of a lifting of the pencil.
For certain families,
they provide a $p$-adic analytic formula for the unit root
as well as the frobenuis polynomial (in the case of rank 4) of each fiber
provided the fiber is ordinary.
Their explicit computations for those examples also reveal
the existence of Dwork type congruences
among the coefficients of the period of such a family.

In this note,
we provide some theoretical explanations
of the observations and computations in their paper.
In particular,
we show that the mod $p$ Dwork congruences
follows from the relative geometry of the family.\medskip

This article is organized as follows.
In \S \ref{DE},
we explore the notion of ordinary differential equations
over a field of characteristic 0 of Calabi-Yau type
and derive some basic algebraic properties of them.
The associated differential modules assemble
the Gauss-Manin connection and the Poincar\'e pairing
on certain parts of the relative de Rham cohomology groups
of families of Calabi-Yau varieties over an affine line
with totally degenerate fibers at the origin.
In \S \ref{Degenerate},
we investigate some cohomological implications
of the existence of a degenerate point
of a pencil of Calabi-Yau varieties.
For example,
we see that the Galois representation
on the cohomology of a certain degenerate fiber
of a family over $\QQ$ coincides with the representation
attached to a modular form.
\S \ref{Hasse} is devoted to the study
of the mod $p$ and $p$-adic properties
of the Calabi-Yau equations.
In particular,
we study the relation between the solutions of a Calabi-Yau
differential equation and the Hasse invariant
of the underlining family modulo $p$.
This relation provides an explanation of the mod $p$ Dwork congruences.
Unfortunately at this stage,
we do not know how to deal with these Dwork congruences in general.
Finally,
some examples are given in \S \ref{Examples},
including a further inspection
of the Hadamard products appeared in \cite{SvS}.

\section{Ordinary differential equations of Calabi-Yau type}
\label{OCYDE}

This section consists of exercises in differential modules.
For basic references, see \cite{vdPS}, Chap 2.
We recall the notion of differential equations of Calabi-Yau type
and derive some basic properties of them.
Such equations arise in the study of the Picard-Fuchs equations
of pencils of Calabi-Yau varieties with certain degenerations.
We postpone the geometric picture to the following sections.
Notice that some properties we derived here
(e.g., Lemma \ref{parity})
can be obtained more easily if the differential equations
are from geometry.\medskip

Let $K$ be a field of characteristic 0
with a fixed embedding into $\CC$.
Denote by $\bar{K} \rc \CC$ the algebraic closure of $K$.
Let $t$ denote a variable.
It will also be regarded as a fixed parameter
of the projective line $\PP^1$ over $K$.
Let $\tha = t\frac{d}{dt}$
be the usual logarithmic differential with respect to $t$.
In this paper,
we adapt the point of view
that the derivative $\tha$ near a regular singular point
is more natural
with respect to the logarithmic structure
associated to the divisor $\{0\} \rc \PP^1$.
For this reason we will use the highly {\em non-standard} convention:
\[ g' := \tha g \]
for a differentiable function $g = g(t)$ of $t$ throughout the discussion.
Now consider an ordinary linear differential operator
of order $n \geq 1$
of the form
\begin{equation}\label{DE}
\LL = \tha^n + \sum_{i=0}^{n-1}a_i \tha^i
\end{equation}
with coefficients $a_i  = a_i(t) \in K(t)$.

\subsection{The condition (N) and the $\beta$-factor}

Suppose $\LL$ is a differential operator of order $n$
of the form (\ref{DE}).
Consider the following condition on $\LL$ :
\begin{description}
\item{(N)} Null exponents:
$\LL$ has a regular singularity at the origin $t = 0$
(i.e., $a_i \in K(t) \dc K\llsb t\rrsb$)
and the associated indicial polynomial of $\LL$
at this point reduces to
\begin{equation}\label{indicial}
s^n + \sum_{i=0}^{n-1} a_i(0) s^i = s^n.
\end{equation}
\end{description}

For any $\LL$ as in (\ref{DE}),
let formally
\begin{equation}\label{beta}
\beta := \exp\left(\frac{2}{n}\int a_{n-1}\frac{dt}{t} \right)
\end{equation}
(i.e., a {\em non-zero} solution of $n\beta' = 2a_{n-1}\beta$).
We call it the {\em $\beta$-factor} of $\LL$.
It exists at least in some differential field extension of $\bar{K}(t)$.

\begin{Lemma}\label{beta_algebraic}
Suppose $\LL$ as in (\ref{DE}) satisfies the condition (N).
Then we have the following.
\begin{enumerate}
\item
The $\beta$-factor in (\ref{beta})
can be taken from $(1+tK\llsb t\rrsb)$.
\item
Regard $\LL$ as a differential operator in $\bar{K}(t)[\tha]$.
Then the differential Galois group of $\LL$
is contained in ${\rm SL}(n)_{\bar{K}}$
if and only if
$\beta$ can be chosen with
\[ \sqrt{\beta^n} \in K(t) \dc \left(1+tK\llsb t\rrsb\right). \]
\end{enumerate}
\end{Lemma}

\pf
By (\ref{indicial}),
$a_{n-1}(0) = 0$
and hence one can formally choose
\[ \int a_{n-1}\frac{dt}{t} \in tK\llsb t\rrsb. \]
Thus the first assertion follows.
The differential Galois group of $\LL \in \bar{K}(t)[\tha]$
is in ${\rm SL_n}$
if and only if
there is a non-zero solution in $\bar{K}(t)$
of the operator $\tha + a_{n-1}$ (\cite{vdPS}, Exercise 1.35.5).
Thus one can require that
\[ \beta \in \sqrt[2/n]{\bar{K}(t)} \dc \left(1+tK\llsb t\rrsb\right) \]
and (ii) follows accordingly.
\qed\\

From now on,
we will always assume that
$\beta \in (1 + tK\llsb t\rrsb)$
if $\LL$ satisfies the condition (N).

\subsection{Self-adjointness and the polarization}

Let $\LL$ be as in (\ref{DE}).
Recall that the {\em formal adjoint} $\LL^*$ of $\LL$
is the differential operator
\[ \LL^* = (-1)^n\tha^n + \sum_{i=0}^{n-1} (-1)^i \tha^i a_i. \]
We say that $\LL$ is {\em self-adjoint}
if, as elements in $K(t)[\tha]$,
\begin{equation}\label{SA}
\LL^* = (-1)^n \beta \LL \beta^{-1},
\end{equation}
where $\beta$ is defined in (\ref{beta}).\medskip

For any $\LL$ of order $n$,
denote by $\MM_\LL$
the left $K(t)[\tha]$-module
with a generator $\eta$ defined by:
\begin{eqnarray}\label{ML}
K(t) &\to&
	\MM_{\LL} := K(t) [\tha]/K(t) [\tha]\LL \\
\nonumber
1 &\mapsto& \eta,
\end{eqnarray}
where the map is the natural projection.
As a $K(t)$-module,
$\MM_{\LL}$ is free of rank $n$
with basis
$\{ \eta^{(i)} \}_{i=0}^{n-1}$,
where
$\eta^{(i)} := \tha^i \eta$.
An element $x \in \MM_{\LL}$ is called {\em horizontal}
if $x' := \tha x = 0$.

Define a filtration $\Fil^{\bullet}$ on $\MM_{\LL}$
by setting $\Fil^i =$ the $K(t)$-submodule generated by
$\{ \eta^{(j)} \}_{j = 0}^{n-1-i}$.
A {\em polarization} on $\MM_{\LL}$ is a
$K(t)$-linear, $(-1)^{n+1}$-symmetric, non-degenerate horizontal pairing
\[ \lan\, ,\ran: \MM_{\LL} \times \MM_{\LL} \to K(t) \]
such that $\lan \Fil^i, \Fil^{n-i} \ran = 0$ for $0 \leq i \leq (n-1)$.
As usual,
we say $\MM_{\LL}$ is {\em polarizable}
if there exists a polarization on it;
$\MM_{\LL}$ is called {\em polarized}
if it is equipped with an underlining polarization.
The aim of this subsection is to prove the following.

\begin{Thm}\label{S=P}
Let $\LL$ be as in (\ref{DE})
and $\beta$ defined in (\ref{beta}).
Then $\LL$ is self-adjoint with $\beta \in K(t)$
if and only if
$\MM_{\LL}$ is polarizable.
\end{Thm}

We first prove the following.

\begin{Lemma}\label{parity}
Let $\lan\, ,\ran$ be a $K(t)$-linear horizontal pairing on $\MM_{\LL}$
such that $\lan \eta, \eta^{(i)} \ran = 0$ for $0 \leq i < (n-1)$.
Let $\gamma = \lan \eta, \eta^{(n-1)} \ran$.
We have
\begin{enumerate}
\item
The pairing is uniquely determined by $\gamma$
and $\gamma = c\beta^{-1}$ for some $c \in K$.
\item
The pairing is $(-1)^{n+1}$-symmetric.
\item
The pairing is a polarization if $\gamma \neq 0$.
\end{enumerate}
\end{Lemma}

\pf
Since $\lan\, ,\ran$ is horizontal,
the values $\lan \eta,\eta^{(i)}\ran, 0 \leq i \leq m$
uniquely determine $\lan \eta^{(i)},\eta^{(j)} \ran$
for all $i + j \leq m$
by a simple inductive procedure of taking derivatives.
Thus in our case,
the pairing is uniquely determined by $\gamma$.
We have $\lan \Fil^i, \Fil^{n-i} \ran = 0$ for $0 \leq i \leq (n-1)$
and the pairing is trivial if $\gamma = 0$.

Since $\lan \eta,\eta^{(n-2)}\ran = 0$,
we have
\begin{eqnarray*}
0 &=& \lan \eta,\eta^{(n-2)}\ran' \\
	&=& \lan \eta',\eta^{(n-2)}\ran + \lan \eta,\eta^{(n-1)}\ran,
\end{eqnarray*}
which implies that
$\lan \eta',\eta^{(n-2)}\ran = -\gamma$.
Similarly
one finds, for $0\leq i\leq n-1$, that
\[ \lan \eta^{(i)},\eta^{(n-1-i)}\ran = (-1)^i\gamma \]
by induction.
In particular,
the pairing is non-degenerate if $\gamma \neq 0$.

Assume $n = 2l$ is even.
Then
\begin{eqnarray*}
\gamma &=& \lan \eta, \eta^{(n-1)} \ran \\
	&=& - \lan \eta', \eta^{(n-2)} \ran \\
	&& \vdots \\
	&=& (-1)^{l-1} \lan \eta^{(l-1)}, \eta^{(l)} \ran.
\end{eqnarray*}
Taking derivatives of the $l$ equations above,
rewriting $\eta^{(n)}$ in terms of $\{\eta^{(i)}\}_{i=0}^{n-1}$ via $\LL$
and summing them up,
one gets
\[ l\gamma' = -a_{n-1}\gamma. \]
Thus by definition,
$\gamma = c\beta^{-1}$ for some $c \in K$.
The computation is similar for $n$ odd.

We now prove the parity of $\lan\, ,\ran$.
We may assume that the pairing is non-degenerate.
Let $\NN = \Hom_{K(t)}(\MM_{\LL}, K(t))$
and consider the natural pairing
\[ (\, ,): \NN \times \MM_{\LL} \to K(t) \]
given by $(f,m) = f(m)$.
We equip $\NN$ with the differential module structure
such that $(\, ,)$ is horizontal.
Then indeed,
as a $K(t)[\tha]$-module,
$\NN$ is generated by $\xi$ with
\[ (\xi, \eta^{(i)}) = \left\{ \begin{array}{ll}
	1 & \text{if $i = n-1$} \\
	0 & \text{otherwise} \end{array}\right. \]
(\cite{vdPS}, Exercise 2.12.6).
Denote $\xi^{(i)} = \tha^i \xi$.
By a similar computation as above,
one deduces that
\[ (\xi^{(i)}, \eta^{(j)}) = \left\{ \begin{array}{ll}
	(-1)^i & \text{if $i + j = n-1$} \\
	0 & \text{otherwise.} \end{array}\right. \]
Now the non-degenerate paring $\lan\, ,\ran$
induces an isomorphism $f$ between $\MM_{\LL}$ and $\NN$
and it sends $\eta$ to $\gamma \xi$.
Thus under $f$,
we can regard $\{ \xi^{(i)} \}_{i=0}^{n-1}$ as another basis of $\MM_{\LL}$.
Moreover, since $\NN$ and $\MM_{\LL}$ are dual to each other
(indeed, $\NN = \MM_{\LL^*}$),
after switching the roles of $\MM_{\LL}$ and $\NN$,
we see that
\[ \lan \xi^{(i)}, \eta^{(j)} \ran = (-1)^{n-1} \lan \eta^{(i)}, \xi^{(j)} \ran. \]
This completes the proof of the parity of $\lan\, , \ran$.
\qed\\

\ni{\em Proof of Thm \ref{S=P}.}
Suppose $\LL$ is self-adjoint and $\beta \in K(t)$.
Then by the Lemma above,
we can equip $\MM_{\LL}$ with the polarization
determined by setting $\lan \eta, \eta^{(n-1)} \ran = \beta^{-1}$.

On the other hand,
suppose $\MM_{\LL}$ is polarized.
Then $\beta \in K(t)$ by the Lemma above.
Multiplying by a non-zero constant,
we may assume that $\lan \eta, \eta^{(n-1)} \ran = \beta^{-1}$.
Then $\beta \eta$ is dual to $\eta^{(n-1)}$
with respect to the basis $\{ \eta^{(i)} \}_{i=0}^{n-1}$.
Thus $\LL^*\beta\eta = 0$ (\cite{vdPS}, Exercise 2.12.6)
and hence $\LL^* = (-1)^n\beta\LL\beta^{-1}$.
\qed

\subsection{Calabi-Yau differential equations and the $q$-coordinate}

\ni{\em Definition.}
A differential operator $\LL \in K(t)[\tha]$
of the form (\ref{DE}) is called {\em Calabi-Yau}
if $\LL$ is self-adjoint
and satisfies condition (N) in \S\ref{OCYDE}($a$).\\

\ni{\em Remark.}
It might be better to called such an $\LL$ as above
{\em locally} or {\em quasi-}Calabi-Yau
since in literature (e.g., in \cite{AZ}, \cite{SvS}),
there are some integral conditions on solutions of $\LL$
(cf.~Thm \ref{F_integral})
and here we ignore the singular types outside the origin.
However there seems no unified definition yet.\\

For any $\LL$ as in (\ref{DE}) satisfying the condition (N),
we set
\begin{equation}\label{Lambda}
\Lam_{\LL} = K\llsb t\rrsb [\tha]/K\llsb t\rrsb [\tha]\LL.
\end{equation}
It is a lattice in the completion
$\MM_{\LL}\tensor_{K[t]} K\llsb t\rrsb$ of $\MM_{\LL}$.
We abuse the notation
by denoting $\Fil^i \rc \Lambda_{\LL}$
the induced filtration of $\Fil^i \rc \MM_{\LL}$.

If furthermore,
$\LL$ is Calabi-Yau.
Denote by $F(t) = F_0(t)$ the unique formal power series solution
of $\LL$ with constant term 1 near $t=0$
\begin{equation}\label{F}
\LL F(t) = 0,\quad F(t) \in 1+ tK\llsb t\rrsb;
\end{equation}
for $n \geq 2$,
denote by $F_i(t), 1\leq i \leq n-1$,
the solutions with logarithmic pole of the form
\begin{equation}\label{G}
\LL F_i(t) = 0,\quad
\sum_{r=0}^{i} (-1)^r F_{i-r}(t)\cdot \frac{\log^r t}{r!} \in K\llsb t\rrsb.
\end{equation}
Thus $\{ F_i \}$ form the Frobenius basis of solutions near $t=0$ of $\LL$.
Let
\begin{equation}\label{wr}
wr_i = wr(F, F_1, \cdots, F_i)
	:= \det \left( F_{r}^{(s)} \right)_{0 \leq r,s \leq i}
\end{equation}
be the wronskians of $\{ F_i\}_{i=0}^i$
and set
\begin{equation}\label{q}
q = \exp\left(\frac{F_1}{F}\right) \in t + t^2K\llsb t\rrsb.
\end{equation}
Thus $K\llsb t\rrsb = K\llsb q\rrsb$.
We call $q$ the {\em $q$-coordinate} of $\LL$.

\begin{Thm}\label{basis_Lambda}
With the notations above,
there exists a unique increasing filtration
$U_i$ of $K\llsb t\rrsb [\tha]$-submodule of $\Lambda_{\LL}$
such that, for all $i$,
\begin{equation}\label{opposite}
\Lambda_{\LL} = U_i \oplus \Fil^{i+1}
\end{equation}
and $U_{i+1}/U_i$ are trivial $K\llsb t\rrsb [\tha]$-modules.
Moreover, up to a constant multiplication,
there exists a unique sequence $\{ u_i \}_{i=0}^{n-1}$
with the following two properties:
\begin{enumerate}
\item
As a $K\llsb t\rrsb$-module,
$U_i$ is generated by $\{ u_r\}_{r=0}^i$.
\item
$u_i$ is of the form $u_i = \sum_{r=0}^{n-1-i} v_{i,r} \eta^{(r)}$
with
\begin{equation}\label{coeff_u}
v_{i,n-1-i} = \beta \frac{wr_i}{wr_{i-1}}.
\end{equation}
\end{enumerate}
Consequently,
we have $u_{i+1}' = \tau_{i+1} \cdot u_i$,
where
\begin{equation}\label{tau}
\tau_{i+1} = \frac{wr_{i-1}wr_{i+1}}{wr_i^2},
\end{equation}
and in particular,
$\tau_1 = (\log q)'$.
\end{Thm}

\pf
We find $\{u_i\}$ satisfying (\ref{opposite}) and (\ref{coeff_u})
by induction.

Since $\LL$ satisfies the condition (N),
up to a constant multiple,
there is a unique non-zero element
\[ u_0 = \sum_{r=0}^{n-1} v_{0,r}\eta^{(r)},\ v_{0,r} \in K\llsb t\rrsb \]
which is horizontal.
Since
\[ \LL\lan u_0, \eta\ran = \lan u_0, \LL\eta\ran = 0, \]
we see that $v_{0,n-1} = c\beta F$ for some constant $c$
by Lemma \ref{parity}.
It is obvious that $c \neq 0$.
Thus after modifying by a constant,
$v_{0,n-1}$ is of the form (\ref{coeff_u}).

Let $\LL^{[i]}$ be the $i$-th exterior product of $\LL$.
Then $\LL^{[i]}$ satisfies the condition (N)
and $\Lambda_{\LL^{[i]}}$
is a quotient of the $i$-th exterior power of
the $K\llsb t\rrsb$-module $\Lambda_{\LL}$.
There exists a unique (up to a scalar) horizontal
$u^{[i]} \in \Lambda_{\LL^{[i]}}$ of the form
\[ u^{[i]} = \sum_{r_1 \leq \cdots \leq r_i} v^{[i]}_{r_1 \cdots r_i}
	\eta^{(r_1)}\wedge\cdots\wedge\eta^{(r_i)} \]
with $v^{[i]}_{n-i,\cdots,n-1} \neq 0$.
(Notice that
$\eta^{(n-i)}\wedge\cdots\wedge\eta^{(n-1)} \neq 0$
in $\Lambda_{\LL^{[i]}}$.)
With the induced pairing,
\[ \LL^{[i]}\lan u^{[i]}, \eta\wedge\cdots\wedge\eta^{(i-1)}\ran =
	\lan u^{[i]}, \LL^{[i]}\left(
		\eta\wedge\cdots\wedge\eta^{(i-1)}\right)\ran = 0. \]
Thus we have, after modifying by a scalar,
$v^{[i]}_{n-i,\cdots,n-1} = \beta^i wr_i \in K\llsb t \rrsb^{\times}$.
Therefore, by subtracting an element in $U_{i-1}$,
we obtain $u_i$ satisfying (\ref{opposite}) and (\ref{coeff_u}).
\qed

\subsection{Examples: lower order cases}

In the remaining of this section,
we consider Calabi-Yau $\LL$ in the form (\ref{DE})
of lower orders more explicitly.\medskip

Suppose $n = 1$.
Then condition (\ref{SA}) is empty.
$\Lam_{\LL}$ is generated by $\beta F\eta$,
which is horizontal.\medskip

Suppose $n = 2$.
Then condition (\ref{SA}) is empty.
$\Lam_{\LL}$ is generated by $\{ u_0, u_1\}$ given by
\begin{eqnarray*}
u_0 &=& \beta[F\eta' - F'\eta], \\
u_1 &=& \frac{\eta}{F}.
\end{eqnarray*}
One computes easily that $u_0' = 0$ and $u_1' = (\log q)' u_0$.
See \S \ref{Examples}$(a)$
for a concrete example.\medskip

Suppose $n = 3$.
We have the following.

\begin{Prop}
A differential operator
$\JJ = \tha^3 + \sum_{i=0}^2 b_i\tha^2 \in K(t)[\tha]$
of order 3 is Calabi-Yau
if and only if
it is the symmetric square
of a Calabi-Yau $\LL$ in the form (\ref{DE})
of order 2.
\end{Prop}

\pf
The operator $\JJ$ is the symmetric square of $\LL$
if and only if
\begin{eqnarray}
\nonumber b_2 &=& 3a_1 \\
\label{3-2} b_1 &=& 4a_0 + a_1' + 2a_1^2 \\
\nonumber b_0 &=& a_0' + 2a_0 a_1.
\end{eqnarray}
From the first two relations in (\ref{3-2}),
we see that
the pair $(b_2,b_1)$ determines the pair $(a_1,a_0)$
uniquely.
On the other hand,
if we rewrite the last equation in (\ref{3-2}) in terms of $b_i$,
then it is equivalent to the condition (\ref{SA}) on $\JJ$.
(Explicitly,
the condition is equivalent to the relation
\begin{equation*}
2\beta b_0 = (\beta b_1)'
	- (\beta b_2)'' + \beta^{(3)},
\end{equation*}
where $\beta$ is the $\beta$-factor (\ref{beta}) of $\JJ$.)
Finally if $\JJ$ is indeed the symmetric square of $\LL$,
it is easy to check that $\JJ$ satisfies (N)
if and only if
$\LL$ does so.
\qed

\begin{Cor}
Let $\JJ = \tha^3 + \sum_{i=0}^2 b_i\tha^2 \in K(t)[\tha]$
be Calabi-Yau of order 3.
There exists a basis $\{ v_i \}_{i=0}^2$ of $\Lam_{\JJ}$
satisfying condition (\ref{opposite}) and (\ref{coeff_u})
in Thm \ref{basis_Lambda}
with
\[ v_2' = (\log \ch{q})' v_1, \]
where $\ch{q}$ is the $q$-coordinate (\ref{q}) of $\JJ$.
\end{Cor}

\pf
By the lemma above,
there is a Calabi-Yau $\LL$ such that
$\Lam_{\JJ}$ is the symmetric square of $\Lam_{\LL}$
as $K\llsb t\rrsb [\tha]$-modules.
Let $\{ u_0, u_1\}$ be the basis of $\Lam_{\LL}$
constructed in the discussion of the case $n=2$ above.
Let
\[ v_0 = u_0^2,\ v_1 = u_0u_1,\ v_2 = \frac{1}{2}u_1^2. \]
Then they form a basis of $\Lam_{\JJ}$
and satisfy the derivative conditions.
\qed\\

\ni{\em Remark.}
In the case of the corollary above,
$\ch{q}$ coincides with the $q$-coordinate of the $\LL$ in the proof.

\subsection{The case $n = 4$ and $5$}

Let us now consider the case when $n = 4$.
Let $\LL$ as in (\ref{DE}) be Calabi-Yau.
Explicitly,
the condition (\ref{SA}) translates to
the relation on the coefficients of $\LL$:
\begin{equation}\label{symp}
a_1 = a_2' + \frac{1}{2} a_2a_3 - \frac{1}{2} a_3''
	- \frac{3}{4} a_3a_3' - \frac{1}{8} a_3^3.
\end{equation}
Via the $\beta$-factor (\ref{beta}) of $\LL$,
equation (\ref{symp}) is equivalent to
\[ \beta a_1 - (\beta a_2)' + (\beta a_3)'' - \beta''' = 0. \]

\begin{Prop}\label{basis}
With the notation above, assume that $\LL$ in (\ref{DE})
is a Calabi-Yau differential equation of order 4.
Let $\beta, F, G:=F_1$ be as given in (\ref{beta}), (\ref{F}), (\ref{G}),
respectively.
Consider the following elements
in the $K\llsb t\rrsb [\tha]$-module $\Lam_\LL$:
\begin{eqnarray*}
u_0 &=& \beta [F\eta''' - F'\eta'' + F''\eta' - F'''\eta]
		+ \beta' [F\eta'' - F''\eta]
		+ (\beta a_2 - \beta'')[F\eta' - F'\eta], \\
u_1 &=& \frac{\beta}{F} \Big[ (FG'-F'G)\eta''
		- (FG''-F''G)\eta' + (F'G''-F''G')\eta\Big], \\
u_2 &=& \frac{F\eta' - F'\eta}{FG'-F'G}, \\
u_3 &=& \frac{\eta}{F}.
\end{eqnarray*}
Then $\{ u_i\}_{i=0}^3$ forms a basis of $\Lam_{\LL}$
and satisfies condition (\ref{opposite}) and (\ref{coeff_u})
in Thm \ref{basis_Lambda}.
We have
\begin{eqnarray*}
u_2' &=& \kappa \cdot (\log q)' \cdot u_1 \\
u_3' &=& (\log q)' \cdot u_2,
\end{eqnarray*}
where $q$ is the $q$-coordinate of $\LL$ and
\[ \kappa = \left(q\frac{d}{dq}\right)^2\left(\frac{F_2}{F}\right)
		\in K\llsb t\rrsb = K\llsb q\rrsb. \]
\end{Prop}

\pf
By direct computation,
we have $u_0' = 0$
and $u_3' = (\log q)' \cdot u_2$.
Also one derives easily that
\begin{eqnarray*}
u_2' &=& \frac{F^2}{\beta(FG'-F'G)^2} \cdot u_1 \\
	&=& \kappa \cdot (\log q)' \cdot u_1,
\end{eqnarray*}
where the second equality comes from \cite{AZ}, Prop 1.
Finally $u_1 = (\log q)' \cdot u_0$ by Thm \ref{basis_Lambda}.
\qed\\

\ni{\em Remark.}
In terms of the $\alpha$-factor introduced in \cite{Y_CY}, Lemma 4.1,
we have (up to a constant multiple)
\[ \alpha = t^3 \beta. \]
\\

We continue to assume that $\LL \in K(t)[\tha]$
with leading coefficient 1 is Calabi-Yau of order 4.
Let
\[ \ch\LL = \tha^5 + \sum_{i=0}^4 \ch{a}_i\tha^i \in K(t)[\tha] \]
be the second exterior power of $\LL$.
Notice that by the self-adjointness of $\LL$,
the operator $\ch\LL$ is of order 5
(see \cite{AZ}, Prop 2 and 3).
Explicitly, we have
\begin{eqnarray}
\nonumber \ch{a}_4 &=& \frac{5}{2}a_3 \\
\nonumber \ch{a}_3 &=& 2a_2 + 2a_3' + \frac{7}{4}a_3^2 \\
\label{b-check} \ch{a}_2 &=& -a_1 + 4a_2' + \frac{7}{2}a_2 a_3 \\
\nonumber \ch{a}_1 &=& -4a_0 + 2a_1' + a_2^2 + a_2''
	+ \frac{3}{2}a_1 a_3 + \frac{3}{2}a_2' a_3 + \frac{1}{4}a_2 a_3^2 \\
\nonumber \ch{a}_0 &=& -2a_0' + a_1'' - 2a_0 a_3 + a_1 a_2
	+ \frac{3}{2} a_1' a_3 + \frac{1}{4} a_1 a_3^2.
\end{eqnarray}

\begin{Prop}
Keep the assumptions and notations as above.
Then $\ch{\LL}$ is Calabi-Yau
and there exists a basis $\{ v_i \}_{i=0}^4$ of $\Lam_{\ch\LL}$
satisfying condition (\ref{opposite}) and (\ref{coeff_u})
in Thm \ref{basis_Lambda} with
\begin{eqnarray*}
v_1' = \kappa \cdot (\log q)' \cdot v_0 &,&
v_2' = (\log q)' \cdot v_1, \\
v_3' = (\log q)' \cdot v_2, &,&
v_4' = \kappa \cdot (\log q)' \cdot v_3,
\end{eqnarray*}
where $q, \kappa$ are defined in Thm \ref{basis}.
Consequently
the $q$-coordinate $\ch{q}$ of $\ch{\LL}$ satisfies
\[ d\log \ch{q} = \kappa\cdot d\log q. \]
\end{Prop}

\pf
The self-adjointness is proved by direct computation
(see below for the explicit relations).
Let $\{ u_i \}_{i=0}^3$ be the basis of $\Lam_{\LL}$
constructed in Thm \ref{basis}.
Being the exterior power of $\Lam_{\LL}$,
the $K\llsb t\rrsb [\tha]$-module $\Lam_{\ch\LL}$
is isomorphic to the quotient
of $\bigwedge_{K\llsb t\rrsb}^2 \Lam_{\LL}$
modulo the condition that
$(u_0 \wedge u_3 - u_1 \wedge u_2)$ is horizontal.
(cf.~\cite{AZ}, Prop 2).
Put
\[ v_0 = u_0 \wedge u_1,\
v_1 = u_0 \wedge u_2,\
v_2 = \frac{1}{2}(u_0 \wedge u_3 + u_1 \wedge u_2),\
v_3 = \frac{1}{2}u_1 \wedge u_3,\
v_4 = \frac{1}{2}u_2 \wedge u_3. \]
One then checks readily that they do the jobs.

On the other hand,
with the notations in \S \ref{OCYDE}($c$),
we have
\begin{eqnarray*}
\log \ch{q} &=& \frac{FF_2' - F'F_2}{FF_1'-F'F_1} \\
	&=& q\frac{d}{dq}\left(\frac{F_2}{F}\right) \\
	&=& \int\kappa\, d\log q.
\end{eqnarray*}
Here we take the integral
congruent to $\log q$ modulo $tK\llsb t\rrsb$.
\qed\\

Conversely,
let $\JJ = \tha^5 + \sum_{i=0}^4 b_i\tha^i \in K(t)[\tha]$
be a Calabi-Yau differential operator of order 5.
If $\beta$ is the $\beta$-factor of $\JJ$,
the self-adjointness of $\JJ$ is equivalent to the following two relations:
\begin{equation}\label{orth}
2(\beta b_2) - 3(\beta b_3)'
	+ 4(\beta b_4)'' - 5(\beta)^{(3)} = 0
\end{equation}
and
\[ 2(\beta b_0) = (\beta b_1)' - (\beta b_2)''
	+ (\beta b_3)^{(3)} - (\beta b_4)^{(4)}
	+ (\beta)^{(5)}, \]
which together are equivalent to the two relations in \cite{SvS}, Prop 2.3.

\begin{Prop}
Let $\JJ$ be Calabi-Yau of order 5 as above.
Then there exists a unique Calabi-Yau $\LL \in K(t)[\tha]$
of the form (\ref{DE}) of order 4
such that $\JJ$ is the second exterior power of $\LL$.
\end{Prop}

\pf
By the first four formulas in (\ref{b-check}),
we can determine the coefficients of $\LL$ uniquely from those of $\JJ$.
One checks that
the condition (\ref{symp}) follows from (\ref{orth}).
The condition (N) on $\LL$ is obvious.
\qed

\section{Degenerations}\label{Degenerate}

For the definitions and basic properties
of logarithmic structures, see \cite{Kato} or \cite{HK}.\medskip

Fix a base field $k$.
Consider a flat projective pencil $\pi: X \to \PP^1$
whose generic fiber is smooth.
We further assume that each singular fiber of $\pi$ is
a union of reduced divisors with normal crossings.
(Over characteristic zero,
this is possible by resolution of singularities
and by passing to a finite steps of base change
of cyclic covering
from $\PP^1$ to $\PP^1$.)
For such a pencil $\pi$,
we equip $X$ and $\PP^1$
with the natural logarithmic structures
associated to the union of singular fibers
(which is a reduced normal crossing divisor on $X$)
and the critical values
(which form a reduced divisor on $\PP^1)$,
respectively.
Then $\pi$ is log-smooth.
Denote by $\om^i = \om^i_{X/\PP^1}$
the (locally free) sheaf on $X$ of relative differential $i$-forms
with log poles with respect to the log structures.

Now suppose that the generic fiber of $\pi$
is an absolutely irreducible Calabi-Yau variety of dimension $m \geq 1$.
We will call such a $\pi$
a {\em nice pencil of Calabi-Yau varieties of dimension $m$}.
Then the sheaf $\pi_*\om^m$ is an invertible sheaf on $\PP^1$.
Suppose there exists a locally direct factor
$\MM$ of $R^m \pi_*\om^{\bullet}$ of rank $(m+1)$
which contains $\pi_*\om^m$ and
is stable under the Gauss-Manin connection $\nabla$.
Now suppose $k = \CC$.
Let $a \in \PP^1(\CC)$ be a $\CC$-valued point
and let $N$ denote
(the logarithm of) the local monodromy around $a$.
Then $N$ is nilpotent and it acts on $\MM$.

We make the following working definition,
which is a special variant of being Hodge-Tate
in the sense of Deligne (\cite{Del_inf}, \S 6).\\

\ni{\em Definition.}
With notations and assumptions as above,
we call $\MM$ {\em totally degenerate} at $a$
if $N^m \neq 0$ on $\MM$.
It is called {\em of rigid type} at $a$
if $N^{m-1} = 0$ but $N^{m-2} \neq 0$ on $\MM$.
We will abuse the notation
by saying that the fiber at $a$ of the family $\pi$
is totally degenerate (resp.~of rigid type)
if there exists an $\MM$ as above
which is totally degenerate (resp.~of rigid type)
at $a$.
We call $\MM$ the {\em degenerate factor} of $\pi$ in this case.\\

\ni{\em Remark.}
Suppose there is a totally degenerate fiber of $\pi: X \to \PP^1$ over $\CC$.
Then the degenerate factor $\MM$
is the unique irreducible locally direct factor
of $R^m\pi_*\om^{\bullet}$
which contains $\pi_*\om^m$
and is stable under $\nabla$.\\

For example, consider the case when $m = 3$.
Let $\MM \lc R^3\pi_*\om^{\bullet}$ be of rank 4,
which is locally a direct summand
and is stable under $\nabla$.
Let $N$ be the local monodromy around a point $a \in \PP^1(\CC)$.
Since
$N^i: \Gr^W_{3+i} \MM_a \to \Gr^W_{3-i} \MM_a(-i)$
is an isomorphism (of Hodge structures)
and the Hodge filtration $\Fil^i$ is locally free for each integer $i$,
there are three possibilities of the degeneration types of the Hodge structure
on $\MM$
at this point:
\begin{enumerate}
\item
No degeneration ($N = 0$ on $\MM$).
\item
The fiber is of rigid type.
In this case,
$h^{1,1} = h^{3,0} = h^{0,3} = h^{2,2} = 1$.
That is,
the Hodge structure $\MM_a$ is a consecutive extension
of the Tate $\CC(-2)$ of weight 4
by a rigid (= rank 2) Calabi-Yau piece of weight 3
by the Tate $\CC(-1)$ of weight 2.
\item
The fiber is totally degenerate.
In this case,
$h^{0,0} = h^{1,1} = h^{2,2} = h^{3,3} = 1$.
That is,
the Hodge structure $\MM_a$ is a consecutive extension
of Tate $\CC(-i)$ of weights $2i = 6, 4, 2, 0$.
\end{enumerate}

\begin{Lemma}\label{non-deg}
Let $\pi: X \to \PP^1$ over $\CC$ be a nice pencil of Calabi-Yau varieties
of dimension $m$ as above with a totally degenerate fiber at $0$.
Then the Poincar\'e pairing is non-degenerate around $0$
on the degenerate factor $\MM$.
\end{Lemma}

\pf
Let $\eta$ be a local basis of sections of $\pi_*\om^m$ near 0.
By assumption, we have
\[ \eta^{(i)} \in \Fil^{m-i} \setminus \Fil^{m-i+1} \]
and they form a local basis of $\MM$ at 0.
Here $\eta^{(i)} = (\nabla(\theta))^i\eta$
($= N^i\eta$ at 0).
Since $(\Fil^i)^{\perp} = \Fil^{m+1-i}$,
the cup-product $\gamma:=\lan \eta, \eta^{(m)}\ran$
is an invertible function near 0.
By Lemma \ref{parity},
the assertion follows.
\qed

\begin{Cor}
Let $\pi: X \to \PP^1$ over $K \rc \CC$
be a nice pencil of Calabi-Yau varieties of dimension $m$
with a totally degenerate fiber at $0$.
Let $\eta$ be a local basis of sections of $\pi_*\om^m$ at $0$
and let $\LL$ be the Picard-Fuchs operator of $\eta$.
Then $\LL$ is a Calabi-Yau differential equation
with respect to the parameter $t$ of order $(m+1)$.
\end{Cor}

\pf
The self-adjointness of $\LL$ follows
from Lemma \ref{non-deg}.
The validity of condition (N) is obvious.
\qed\\

Now suppose that the base field $k = \QQ$ and
the nice pencil $\pi: X \to \PP^1$ over $\QQ$
of Calabi-Yau varieties of dimension $m$
is of rigid type at $a \in \PP^1(\QQ)$.
Denote by $W_{\bullet}$
the corresponding monodromy filtration
of the degenerate factor $\MM_a$ at $a$.
Then
\[ \dim_{\QQ} W_{i}/W_{i-1} = \left\{ \begin{array}{cll}
	0 && \text{if $i$ odd and $i \ne m$} \\
	1 && \text{if $i$ even, $2\le i \le 2m-2$ and $i \ne m$} \\
	2 && \text{if $i = m$ is odd} \\
	3 && \text{if $i = m$ is even}. \end{array}\right.\]
If $m \geq 3$,
the subquotient $\NN = W_m/(W_{m-1} + NW_{m+2})$
is then of rank two.

\begin{Prop}
The ${\rm Gal}(\bar{\QQ}/\QQ)$-representation
on the \'etale realization $\NN_{et}$ of $\NN$ coincides
with the representation attached to a cusp form
of weight $(m+1)$.
\end{Prop}

\pf
One knows $\NN_{et}$ is irreducible
by \cite{Se}, (4.8.9).
If $m$ is odd, the statement is a consequence of Serre's conjecture
(\cite{Se}, ($3.2.4_?$) and Th 6),
which has recently been proved in \cite{KW}.
If $m$ is even, one uses \cite{L}, Cor 1.4.
\qed

\section{Mod $p$ and $p$-adic aspects}\label{Hasse}

In this section,
we fix a prime $p$ and suppose $p > m$.
Let (for simplicity) $K$ be a finite unramified extension of $\QQ_p$
with ring of integers $W$ and residue field $k$.
Let $\pi:X \to \PP^1$ be a nice pencil of Calabi-Yau varieties
of dimension $m$ over $K$
with a totally degenerate fiber at 0.
Let $\om^i$ be the sheaf of relative differential $i$-forms
with log-poles along the logarithmic structure.
We fix a basis $\eta$ near $0$ of sections of
$\pi_*\om^m \rc R^m\pi_*\om^{\bullet}$
and let $\LL$ be the Picard-Fuch operator of $\eta$
as before.
We call that $\pi$ has {\em nice} reduction
if $\pi$ has a flat model over $W$
such that the reduction $\bar{\pi}: \bar{X} \to \PP^1$ over $k$
is also a nice pencil of Calabi-Yau varieties.

\subsection{The $p$-adic input}

\begin{Lemma}\label{ord}
Suppose the pencil $\pi$ has nice reduction
and the degenerate factor $\MM_0$ at $0$ is stable
under (a lift of) the absolute frobenius.
Then the frobenius action on $\MM_0$ is ordinary
and consequently $\MM$ is generically ordinary.
\end{Lemma}

\pf (Cf.~\cite{Y_Dw}, Thm 2.2.)
Since $\MM$ is totally degenerate at 0,
the Hodge structure on $\MM_0$
is a successive extension of rank 1 Hodge structures
of pairwise different weights.
Since the frobenius respects the Hodge and the weight filtrations
(cf.~\cite{Mo}, Remarques 3.28),
the result follows.
\qed

\begin{Prop}
With the assumptions and notations in Lemma \ref{ord}
and $\tau_i$ as in (\ref{tau}),
we have $\tau_i = \tau_{m+1-i}$ for all $1\leq i\leq m$.
\end{Prop}

\pf
Let $\Lambda_{0,cris}$ be the $m$-th logarithmic crystalline cohomology
of the fiber of $\bar{\pi}$ over the localization at $0$.
Retain the notations in Thm \ref{basis_Lambda}
and let $\phi$ be a lift of the absolute frobenius on $\Lambda_{0,cris}$.
Then $u_i \in \Lambda_{0,cris}$
by (the log version of) \cite{Katz_Dw}, Prop 3.1.
It is also known that
$\lan x,y \ran$ is divisible by $p^{2m}$
for any $x,y \in \Lambda_{0,cris}$
and
\[ \phi(\eta^{(i)}) \in
	p^{m-i}\Lambda_{0,cris} \setminus p^{m-i+1}\Lambda_{0,cris} \]
by the Lemma above.
Thus $\lan u_i, u_{m-j} \ran = 0$ if $i \neq j$
and $\lan u_i, u_{m-i} \ran$ are non-zero constants.
Consequently,
\begin{eqnarray*}
0 &=& \nabla(\tha) \lan u_i, u_{m+1-i} \ran \\
	&=& \tau_i \cdot \lan u_{i-1}, u_{m+1-i} \ran
		+ \tau_{m+1-i} \cdot \lan u_i, u_{m-i} \ran.
\end{eqnarray*}
Since $\tau_i(0) = 1$ for all $i$,
the assertion follows.
\qed\\

\ni{\em Question.}
Can one derive the above Proposition
as a property of Calabi-Yau differential equations
without referring the underlying families?\\

\begin{Lemma}\label{betaF}
Suppose the pencil $\pi$ over $K$ with a totally degenerate fiber at $0$
has nice reduction
and the degenerate factor $\MM_0$ at $0$ is stable under frobenius.
Then there exists a non-zero constant $c \in W$ such that
$c\beta F \in W\llsb t\rrsb$.
\end{Lemma}

\pf
By Lemma \ref{ord} and \cite{Katz_Dw}, Prop 3.1.3
(cf.~\cite{Y_CY}, Cor 2.2),
there is a non-zero element
\[ u = \sum_{i=0}^m v_i\eta^{(i)},\ v_i \in W\llsb t\rrsb \]
which is horizontal.
By Thm \ref{basis_Lambda},
$v_m = c\beta F$ for some constant $c$.
\qed

\begin{Thm}\label{F_integral}
Let $\pi$ be a nice pencil of Calabi-Yau varieties of dimension $m$
over $\QQ$.
Fix a basis near $0$ of sections of $\pi_*\om^m$
and let $\LL$ be the associated Picard-Fuchs operator.
Let $F$ be the formal solution of $\LL$ as in (\ref{F}).
Then $F(t) \in (1+t\ZZ_p\llsb t\rrsb)$ for all $p$ sufficiently large.
In particular,
$F(t) \in \ZZ[1/N]\llsb t\rrsb$ for some integer $N \neq 0$.
\end{Thm}

\pf (Cf.~\cite{SvS}, Conj 2.1.)
By the Lemma above and its proof,
we see that except for a finite set of primes,
$\beta F \in \ZZ_p\llsb t\rrsb$
for all $p$ such that the pencil $\pi$ has nice reduction.
By Lemma \ref{parity},
$\beta(t)$ is rational over $\QQ$
and hence is in $(1 + t\ZZ_p\llsb t\rrsb)$
for almost all $p$.
Thus the assertion follows.
\qed

\subsection{The Hasse invariant and the differential equation}

In the remaining of this section,
we keep the assumptions and notations in the beginning of \S \ref{Hasse}
and assume that $\pi$ has nice reduction $\bar{X}$.
Let $\bar\om^i$ be the sheaf on $\bar{X}$
of relative differential $i$-forms with log-poles.
We choose $\eta$ that can be extended to the flat model of $\pi$
and assume the following condition is satisfied:
\begin{description}
\item{(R)}
Non-degeneracy of the reduction of $\eta$ at the origin:
\[ \lan \eta, \eta^{(m)}\ran \in W^{\times} + tW\llsb t\rrsb. \]
\end{description}
Let $\sigma$ be the lift of the absolute frobenius on $W$
to $W\llsb t\rrsb$ by sending $t$ to $t^p$.
For any $x \in W\llsb t\rrsb$,
denote by $x^{\sigma}$ the image under $\sigma$.\medskip

Let $\bar{\HH} = \bar{\pi}_*\bar\om^m$.
Consider the adjoint morphism
\[ V: \bar\HH \to \sigma^*\bar\HH \]
of the absolute frobenius with respect to the cup-product
on $R^m\pi_*\bar\om^{\bullet}$.
Represent $V$ by $\Has$
defined by
\[ V(\eta) = \Has \cdot \sigma^*\eta. \]
It can be regarded as an element in $k(t)$.
We call $\Has$ the {\em Hasse invariant} of the family $\bar{\pi}$
(with respect to $\eta$).
Notice that $\Has(0)$ is well-deined and non-zero at $t = 0$
by Lemma \ref{ord}.

\begin{Prop}
We have $\LL \Has = 0$.
That is,
with respect to the parameter $t$,
the Hasse invariant $\Has$ is a rational solution of $\LL$
in characteristic $p$.
\end{Prop}

\pf
Since $V$ is horizontal, $\nabla(\LL)\eta = 0$
and $\nabla \sigma^* \eta = 0$,
we see that $\LL \Has = 0$.
\qed\\

Thus, for example,
if
$\pi_*\bar\om^m = \OO(1)$ on $\PP^1$,
then $\Has$ is a section of $\OO(p-1)$.
Consequently,
if $\eta \in \Gamma(\PP^1, \pi_*\bar\om^m)$
and regarding $\Has$ as a function of $t$,
we see that $\Has \in k[t]$ is of degree $\leq (p-1)$.\medskip

If (R) is satisfied,
then by Lemma \ref{parity}, \ref{betaF} and their proofs,
$F \in W\llsb t\rrsb$.
Let $f(t) = F(t)/F(t)^{\sigma}$
regarded as a formal power series over $W$.
Let $\bar{f}(t)$ be the reduction mod $p$ of $f(t)$.

\begin{Prop}
Suppose $\pi$ satisfies the condition (R) above.
Then we have the following.
\begin{enumerate}
\item
$\bar{f}(t) \in (1 + tk\llsb t\rrsb) \dc k(t)$.
\item
The function $\bar{f}$ satisfies $\LL \bar{f} = 0$.
\end{enumerate}
\end{Prop}

\pf
By \cite{Katz_Dw}, 4.1.9,
the function $\beta F/(\beta F)^{\sigma}$
is a lift of an element in $k(t)$.
Thus the first assertion follows.

Formally we have
\[ F(t) = \frac{F(t)}{F(t)^{\sigma}}
	\frac{F(t)^{\sigma}}{F(t)^{\sigma^2}}\cdots
	= f(t)f(t)^{\sigma}\cdots. \]
Applying the differential operator $\LL$,
we have
\begin{eqnarray*}
0 = \LL F(t) &=& \LL \big(f(t)f(t)^{\sigma}\cdots\big) \\
	&\equiv& \big(\LL f(t)\big)\big(f(t)^{\sigma} \cdots\big) \\
	&\equiv& \big(\LL \bar{f}(t)\big)\big(\bar{f}(t)^p \cdots\big) \pmod{p}.
\end{eqnarray*}
Thus the rational function $\bar{f}(t)$ is a solution of $\LL$
in characteristic $p$.
\qed

\begin{Prop}
Assume condition (R) is fulfilled.
Let $c = \Has(0)$.
We have $\Has = c\bar{f}$
regarded as rational functions of $t$.
\end{Prop}

\pf
Over a non-empty open subset of $\PP^1$,
both of the two functions $\Has$ and $c\bar{f}$
represent the absolute frobenius action on $\eta$
(\cite{Katz_Dw}, 4.1.9).
Thus the assertion follows.
\qed

\begin{Cor}
Assume that condition (R) is fulfilled
and $\pi_*\om^m = \OO(1)$.
Suppose $\eta \in \Gamma(\PP^1, \pi_*\bar\om^m)$
and let $c = \Has(0)$.
Then $c\bar{f}(t) = \Has(t) \equiv cF^{<p}(t) \mod{p}$,
where $F^{<p}(t)$ is the truncation of $F(t)$
up to degree $(p-1)$.
\end{Cor}

\pf
Under the assumptions,
$\Has(t)$ is a polynomial of degree $< p$.
\qed\\

\ni{\em Remark.}
The statement of the Corollary
is equivalent to the mod $p$ case
of the Dwork congruences in \cite{SvS}, \S 2.3.

\subsection{The higher Hasse invariant}

Suppose that $\LL$ satisfies the condition (R)
and for simplicity that $R^m\pi_*\om^{\bullet}$ is of rank $(m+1)$.
Let $\MM_{cris}$ be the
relative $m$-th logarithmic crystalline cohomology of $\bar{\pi}$.
Assume $p$ odd.
Consider $\NN:= \bigwedge^2 \MM_{cris}$
equipped with the induced cup-product pairing
and with frobenius $= p^{-1} \cdot
(\text{the induced frobenius from $\MM_{cris}$
to $\bigwedge^2 \MM_{cris}$})$.
Then $\xi = \eta\wedge\eta'$ is a local section of $\NN$
and $\not\equiv 0 \mod{p}$.\medskip

Let $\bar\NN = (\NN \mod{p})$
and $V$ the adjoint of the frobenius on $\bar\NN$.
Define the Hasse invariant $\ch\Has$
(with respect to $\xi$) of $\NN$ by
\begin{eqnarray*}
V: \bar\NN &\to& \sigma^*\bar\NN \\
	\xi &\mapsto& \ch{\Has}\cdot\sigma^*\xi.
\end{eqnarray*}
If $\eta$ is a section over an open subset $U$ of $\PP^1$,
then we have
\[ \ch\Has \in \Gamma\left(U,
	\left(\pi_*\bar\om^m \tensor R^1\pi_*\bar\om^{m-1}\right)\right), \]
and for $x \in U(\bar{k})$,
the Newton polygon of $\MM_{cris}$ over
$x$ starts with slopes
\[\left\{ \begin{array}{cll}
	0, 1 &&
		\text{if $\Has(x) \neq 0$ and $\ch\Has(x) \neq 0$} \\
	0, >1 &&
		\text{if $\Has(x) \neq 0$ but $\ch\Has(x) = 0$} \\
	1/2, 1/2 &&
		\text{if $\Has(x) = 0$ but $\ch\Has(x) \neq 0$} \\
	> 1/2 &&
		\text{if $\Has(x) = 0$ and $\ch\Has(x) = 0$}. \end{array}\right.\]
Consequently
the variation of crystals $\MM_{cris}$ (over $U$)
away from $\ch\Has = 0$
is an extension (\cite{Katz_Slope}, Thm 2.4.2)
\[ 0 \to \MM_{\leq 1} \to \MM_{cris}\big|_{U \setminus \{\ch\Has = 0\}}
	\to \MM_{> 1} \to 0, \]
where $\MM_{\leq 1}$ (resp.~$\MM_{>1}$)
is the slope $\leq 1$ (resp.~$>1$) part
which is of rank 2.\medskip

On the other hand,
let
\[ \ch{F}(t) = (FG'-F'G)(t) \in K\llsb t\rrsb. \]
Then the condition (R) on $\LL$
implies that indeed, $\ch{F}(t) \in W\llsb t\rrsb$.
Similarly to the discussion in the previous subsection,
if $\pi_*\om^m \tensor R^1\pi_*\om^{m-1} = \OO(1)$
and $\xi$ is a global section,
we have
$\ch\Has(t) \equiv \ch{c}\ch{F}^{<p}(t) \mod{p}$,
where $\ch{c} = \ch\Has(0)$.\\

\ni{\em Question.}
Can one detect the sheaves
$\pi_*\om^m$ and $R^1\pi_*\om^{m-1}$ on $\PP^1$
from the Picard-Fuchs operator?
Is there a geometric/homological interpretation
of the higher Dwork congruences?

\section{Examples}\label{Examples}

\subsection{The Legendre family}

Let $\lam$ be a fixed parameter of $\PP^1$.
Consider the Legendre family $\pi: E \to \PP^1$ over $\QQ$
whose affine part is given by
\[ y^2 = x(x-1)(x-\lam). \]
The Picard-Fuchs operator $\LL$ associated to the invariant differential
\[ \eta = \frac{dx}{2y} \]
is the one associated to the Gauss hypergeometric series
\[ F(\lam) := {}_2F_1(\frac{1}{2}, \frac{1}{2}, 1; \lam). \]
Note that the monodromy around $\lam = \infty$
is not unipotent.
Consider the double cover $a: \PP^1 \to \PP^1$
given by $\lam = t^{-2}$.
Then
\[ 2\eta = \frac{dx}{\sqrt{x(x-1)(x-1/t^2)}}
	= \frac{t dx}{\sqrt{x(x-1)(t^2 x-1)}}, \]
which is zero precisely when $t = 0$.
Thus with respect to $t$,
we have $(a^*\pi)_*\omega^1_{a^*E/\PP^1} = \OO(1)$.\medskip

Let $p$ be an odd prime.
the discussion above shows that
the Hasse invariant of the family over $\FF_p$
in the affine part with respect to $\lam$
is a polynomial of degree $\leq \frac{p-1}{2}$.
In fact,
\[ \Has(\lam) \equiv (-1)^{(p-1)/2}F^{<p}(\lam) \mod{p} \]
is of degree exactly $\frac{p-1}{2}$.
This is due to the fact that
$\Has(\lam)$ has only simple roots (\cite{Si}, Thm V.4.1)
and that the singular curve corresponding to $t = 0$
is ordinary.\medskip

Let an upper $'$ denote the derivative $d/d\log \lam$.
The set
\[ \left\{ u_0 = (1-\lam)[F\eta' - F'\eta],
	u_1 = \frac{\eta}{F} \right\}, \]
provides a local basis of
\[ \HH := R^1\pi_*\left(\om^{\bullet}_{E|_{\AA^1}/\AA^1}\right) \]
near $\lam = 0$
adapted to the slope filtration
(cf.~\cite{Katz_Dw}, \S 8).
\medskip

The Galois representation
on the cohomology $\HH_0$ over $\QQ$ at $\lam = 0$ is reducible.
It has the form
\[ 0 \to \QQ(\epsilon) \to \HH_0 \to \QQ(-1)(\epsilon) \to 0, \]
where $\epsilon$ is the Legendre character.
This is simply because
the corresponding singular curve splits over $\QQ(\sqrt{-1})$
but not over $\QQ$.
This gives an explanation of the constant term of $\Has(\lam)$.

\subsection{Dwork families}

Fix an integer $n \geq 2$.
The Dwork family of Calabi-Yau varieties of dimension $(n-1)$
is the pencil of hypersurfaces in $\PP^n$ given by the equation
\[ {\cal P}_t:
	X_1^{n+1} + \cdots + X_{n+1}^{n+1} - (n+1)t X_1 \cdots X_{n+1}. \]
In this case,
we consider the differential
\[ \eta = {\rm Res}\, \frac{t\Omega}{{\cal P}_t}, \]
where
\[ \Omega = \sum (-1)^i X_i dX_1 \wedge \cdots
	\wedge \widehat{dX_i} \wedge \cdots \wedge dX_{n+1}. \]
Via the parameter $\lam = t^{-(n+1)}$,
the associated Picard-Fuchs operator $\LL$
is the one associated to the generalized hypergeometric series
\[ {}_nF_{n-1}\left(\begin{array}{c}
	\frac{1}{n+1}, \cdots, \frac{n}{n+1}\\
	1, \cdots, 1 \end{array}; \lam \right). \]

Explicitly one can pick the annihilator of $\eta$ to be
\begin{eqnarray*}
\LL &=& \tha_{\lam}^n
		- \lam\prod_{i=1}^n\left( \tha_{\lam} + \frac{i}{n+1} \right) \\
	&=& (1-\lam)\tha_{\lam}^n
		- \frac{n}{2}\lam\tha_{\lam}^{n-1} + \cdots
\end{eqnarray*}
where $\tha_{\lam} = \lam\frac{d}{d\lam}$.
One checks that
$\beta = (1-\lam)$ is the $\beta$-factor of $\LL$
regarded as a Calabi-Yau differential operator
in $\QQ(\lam)[\tha_{\lam}]$.\medskip

Similar to the case of the Legendre family,
the monodromy at $\lam = \infty$ is not unipotent.
However with respect to $t$,
the form $\eta$ is well-defined everywhere
and vanishes precisely when $t = 0$.
Thus $\pi_*\om^{n-1} = \OO(1)$.
Consequently for the reduction over $\FF_p,\, p \nmid (n+1)$,
the Hasse invariant $\Has(\lam)$,
regarded as a polynomial of $\lam$ here,
is a polynomial of degree at most $\lfloor \frac{p-1}{n+1} \rfloor$.
In fact the degree of $\Has(\lam)$ is exactly the upper bound.
This is because the Fermat point $t=0$
is not ordinary if and only if $p \not\equiv 1 \mod{(n+1)}$
and there is no $n$-multiple root of ``$\Has(t)$".
Notice that in this example,
double roots do occur in $\Has(\lam)$.\medskip

On the other hand,
over each geometric point of
$\mu_{n+1} := \Spec \QQ[t]/(t^{n+1}-1)$ in $\PP^1$,
the fiber of the Dwork family has ordinary double points
as its singularities.
Thus $N^2 = 0$ for the local monodromy $N$
around a point of $\mu_{n+1}$.
One can show that $N \ne 0$ in this case
(\cite{HST}, Lemma 1.6).
Now consider the fiber over $t=1$ in the case $n = 4$.
Retain the notations in the end of \S \ref{Degenerate}.
Then the fiber is of rigid type
and the subquotient $W_3/W_2$ is modular,
which was first proved by Schoen in \cite{Sc}.
The corresponding modular form is of weight 4 and level 25
and with the trivial character.

\subsection{Hadamard products}

Here we describe how to obtain the unit roots precisely
for certain Hadamard products considered in \cite{SvS}, \S 3.
The only missing piece in loc.cit.~is
to determine the constant $\epsilon_4$ in Prop 2.7 there.
To do this,
we study the frobenius action
on the cohomology of the totally degenerate fiber
by applying the weight spectral sequence in \cite{Mo}.
For references of Hadamard products
and examples of pencils of elliptic curves we discuss here,
see \cite{SvS}, \S\S 3.1 and 3.2.\medskip

Let $X, Y \to \PP^1$ be two pencils of elliptic curves
over a finite field $k$ of characteristic $p$
with totally degenerate fibers $X_0, Y_0$ at 0, respectively.
We assume that $X_0, Y_0$
are strictly normal crossing divisors.
Let $\xi_1$ and $\xi_2$ be
local bases at 0 of horizontal sections
of the relative $H^1_{cris}$
of $X$ and $Y$ over $\PP^1$, respectively.
Then $\xi_i$ are eigenvectors of the relative frobenius.
Denote by $c_i$
the corresponding eigenvalues.
Then $c_1 = 1$
if the degenerate curve $X_0$ is of split multiplicative type over $k$;
$c_1 = -1$ if $X_0$ is non-split.
The statement for $c_2$ is similar.\medskip

Geometrically,
the Hadamard product comes from the following diagram
with all squares being Cartesian:
\[\xymatrix{
Z \ar[r]\ar[d] &
	\widetilde{X\times Y} \ar[r]\ar[d] & X\times Y \ar[d] \\
C_1 \uc C_2 \ar[r]\ar[d] &
	\widetilde{\PP^1\times \PP^1} \ar[r]^b\ar[d] & \PP^1\times \PP^1 \\
	0 \ar[r] & \PP^1. }\]
Here $b$ is the blow-up of $\PP^1 \times \PP^1$
along $(\infty,0)$ and $(0,\infty)$;
$C_i$ are rational curves
with $C_1 =$ the strict transformation of $b^{-1}(\AA^1 \times 0)$
and $C_2 =$ that of $b^{-1}(0\times \AA^1)$.
Notice that
the map $b$ induces isomorphisms
from $C_1$ to $\PP^1 \times 0$
and from $C_2$ to $0 \times \PP^1$
and $C_1 \dc C_2 = 0:=b^{-1}(0,0) \in \widetilde{\PP^1\times \PP^1}$.
Over $\bar{k}$,
write $X_0 = \bigcup D_i$ and $Y_0 = \bigcup E_j$,
where $D_i$ and $E_j$ are distinct projective lines.
Let $Z^{(i)}$ be the disjoint union of
all possible intersections of $i$ distinct irreducible components of $Z$.
We then have
\begin{eqnarray}
\nonumber
Z^{(1)} &=& \left(\bigsqcup D_i \times Y\right)
		\sqcup\left(\bigsqcup X \times E_j\right); \\
\nonumber
Z^{(2)} &=& \left(\bigsqcup (D_i\dc D_r) \times Y\right)
		\sqcup\left(\bigsqcup D_i \times E_j\right)
		\sqcup\left(\bigsqcup X \times (E_j\dc E_s)\right); \\
\nonumber
Z^{(3)} &=& \left(\bigsqcup (D_i\dc D_r) \times E_j\right)
		\sqcup\left(\bigsqcup D_i \times (E_j\dc E_s)\right); \\
Z^{(4)} &=& \bigsqcup \big((D_i\dc D_r) \times (E_j\dc E_s)\big).
\end{eqnarray}
\medskip

Let us recall the weight spectral sequence
(\cite{Mo}, \S 3.23; cf.~\cite{St}, Cor 4.20):
\[ E_1^{-j,i+j} = \bigoplus_{r\geq 0,r\geq -j}
	H_{cris}^{i-j-2r}(Z^{(1+j+2r)}/W)(-j-r)
	\Longrightarrow H^i(Z^{\times}/W^{\times}), \]
which degenerates at $E_2$ modulo torsion (\cite{Mo}, Th 3.32).
Here $W$ is the ring of Witt vectors of $k$
and the target $H^i(Z^{\times}/W^{\times})$
is the $i$-th logarithmic crystalline cohomology of $Z$.
Now assume that $X$ and $Y$ have trivial crystalline cohomology
of odd degrees
and they are ordinary.
Then the weights of the $E_1$-terms are all integers
and
the non-zero terms of $E_1$ appear only when the weight are even.
Thus,
putting $E_1^{r,s}$ at the $(r,s)$-spot,
the complete picture of the $E_1$-terms looks like
\[ \xymatrix{
\text{wt 6:} & E_1^{-3,6} & E_1^{-2,6} & E_1^{-1,6} & E_1^{0,6} \\
\text{wt 4:} && E_1^{-2,4} & E_1^{-1,4} & E_1^{0,4} & E_1^{1,4} \\
\text{wt 2:} &&& E_1^{-1,2} & E_1^{0,2} & E_1^{1,2} & E_1^{2,2} \\
\text{wt 0:} &&&& E_1^{0,0} & E_1^{1,0} & E_1^{2,0} \ar[r]^d & E_1^{3,0}. }\]

\begin{Lemma}
Let $\xi$ be the Hadamard product of $\xi_i$
and $c$ the eigenvalue of the relative fronbenius action on $\xi$.
Then $c = c_1 c_2$.
\end{Lemma}

\pf
$c$ represents the relative frobenius action on the cokernel of $d$,
where $d$ is the boundary map
in the displayed $E_1$-terms above.
Let $K$ be the field of fractions of $W$.
By \cite{Mo}, Lemme 5.2,
we see that $({\rm coker}\ d)\tensor K$ is 1-dimensional
and the frobenius acts on $\xi$ as the product
of its actions on $\xi_i$.
Thus the statement follows.
\qed

Department of Mathematics

National Taiwan University

Taipei, 10617 Taiwan

E-mail address: {\tt jdyu@math.ntu.edu.tw}
\end{document}